\documentclass[a4paper,12pt]{article}
\usepackage[english]{babel}
\usepackage[T2A]{fontenc}
\usepackage[cp1251]{inputenc}
\usepackage{amsthm}
\usepackage[tbtags]{amsmath}
\usepackage{amsfonts,amssymb}
\begin{document}


\begin{center}
\textbf{E. I. Kompantseva, A. A. Tuganbaev}

\textbf{Multiplication Groups \\ 
of Abelian Torsion-Free Groups of Finite Rank}
\end{center}

\textbf{Abstract.} For an Abelian group $G$, any homomorphism $\mu\colon G\otimes G\rightarrow G$ is called 
a \textsf{multiplication} on $G$. The set $\text{Mult}\,G$ of all multiplications on an Abelian group $G$ itself is an Abelian group with respect to addition; the group is called the \textsf{multiplication group} of $G$.  Let $\mathcal{A}_0$ be the class of all reduced block-rigid almost completely decomposable groups of ring type with cyclic regulator quotient. In this paper,  for groups  $G\in \mathcal{A}_0$, we describe groups $\text{Mult}\,G$. We prove that for $G\in \mathcal{A}_0$, the group $\text{Mult}\,G$ also belongs to the class $\mathcal{A}_0$.
For any group $G\in \mathcal{A}_0$, we describe the rank, the regulator, the regulator index, invariants of near-isomorphism, a main decomposition, and a standard representation of the group $\text{Mult}\,G$.


\textbf{MSC2020 datebase:} 20K30, 20K99, 16B99 

\section{Introduction}

For an Abelian group $G$, a \textsf{multiplication} on $G$ is a homomorphism $\mu\colon G\otimes G\rightarrow G$. The set $\text{Mult}\,G$ of all multiplications on the group $G$ itself is an Abelian group with respect to addition; the group is called the \textsf{multiplication group} of $G$ or the \textsf{group of multiplications} on $G$ \cite{Fuc15}. An Abelian group $G$ with multiplication on $G$ is called a \text{ring on the group} $G$. The problem of studying the relationship between the structure of an Abelian group and the properties of ring structures on it is very multifaceted and has a long history in algebra; see\cite{AndW17}, \cite{BeaP61}, \cite{Fei97}, \cite{Fei00}, \cite{Gar74}, \cite{Jac82}, \cite{Kom10}, \cite{Kom15}. 

In this paper, we consider only additively written Abelian groups and "a group" means "an Abelian group" in what follows.

In this paper, we study the group $\text{Mult}\,G$ for an almost completely decomposable Abelian group $G$. 
A torsion-free group $G$ of finite rank is called an \textsf{almost completely decomposable} group (\textsf{$ACD$-group}) if $G$ contains completely decomposable subgroup of finite index. $ACD$-groups were studied in \cite{BlaM94}, \cite{Bur84}, \cite{DugO93}, \cite{Kom09}, \cite{Lad74}, \cite{Mad00} and other papers. 
The achieved level of $ACD$-group theory development is recorded in the book \cite{Mad00}. 

Any $ACD$-group $G$ contains a special uniquely defined completely decomposable (see \cite{Fuc15}) subgroup $\text{Reg}\,G$ of finite index which is a fully invariant subgroup of $G$; it is called the \textsf{regulator} of the group $G$. The regulator of an $ACD$-group can be defined as the intersection of all its completely decomposable subgroups of lowest index \cite{Bur84}.
The factor group $G/\text{Reg}\,G$ is called the \textsf{regulator quotient} of the group $G$; the index of the subgroup $\text{Reg}\,G$ in the group $G$ is called a \textsf{regulator index}. It is denoted by $n(G)$. 
$ACD$-groups with cyclic regulator quotient are often called \textsf{$CRQ$-groups}.

Let $G$ be an almost completely decomposable group. Then the group $\text{Reg}\,G$ can be uniquely, up to isomorphism, represented as a direct sum of torsion-free groups of rank 1 \cite[Proposition 86.1]{Fuc73}. For every type $\tau$, we denote by $\text{Reg}_{\tau}\,G$ the sum of summands of rank 1 and type $\tau$ in this decomposition of the group $\text{Reg}\,G$. The set of types 
$$
T(G)=T(\text{Reg}\,G)=\{\tau\,|\,\text{Reg}_{\tau}\,G\ne 0\}
$$
is called the \textsf{set of critical types} of groups $G$ and $\text{Reg}\,G$. 
If $T(G)$ consists of pairwise incomparable types, then the groups $G$ and $\text{Reg}\,G$ are called \textsf{block-rigid} groups. If, in addition, for any $\tau\in T(G)$, the group $\text{Reg}_{\tau}G$ is of rank $1$, then $G$ and $\text{Reg}\,G$ are called \textsf{rigid} groups. If all types in $T(G)$ are idempotent types, then $G$ is called a group \textsf{of ring type}.

We note that a block-rigid $ACD$-group is either divisible or reduced. For a divisible torsion-free group $G$, the group $\text{Mult}\,G$ is described in \cite[Section~121]{Fuc73}; therefore, we consider only reduced groups in what follows.

We denote by $\mathcal{A}_0$ the class of all reduced block-rigid $CRQ$-groups of ring type. In Section 2, we describe the group $\text{Mult}\,G$ for $G\in\mathcal{A}_0$ (Theorem 2.8).
The aim of Section 3 is, for groups $G$ in the class $\mathcal{A}_0$, to study properties of $\text{Mult}\,G$. It is proved (Theorem 3.4) that if $G$ is a block-rigid $CRQ$-group of ring type, then $\text{Mult}\,G$ also is a block-rigid $CRQ$-group of ring type. We describe the rank, the regulator, the regulator index, invariants of near-isomorphism, a main decomposition and a standard representation of the group $\text{Mult}\,G$ for $G\in\mathcal{A}_0$.

The multiplication $\mu\colon G\otimes G\rightarrow G$ is often denoted by the symbol $\times$, i.e.,
$$
\mu(g_1\otimes g_2)=g_1\times g_2 \text{ for all } g_1,g_2\in G.
$$
The multiplication $\times$ on the group $G$ induces a ring on this group which is denoted by $(G,\times)$.
Let $G$ be a group and $g\in G$. The characteristic and the order of the element $g$ are denoted by $\chi(g)$ and $o(g)$, respectively. The rank and the divisible hull of the group $G$  are denoted by $r(G)$ and $\tilde{G}$, respectively. 
If $S\subseteq G$, then $|S|$ is the cardinality of the set $S$ and $\langle S\rangle$ is the subgroup of the group $G$ generated by the set $S$. We write an element of a group direct product $\prod_{i\in I}G_i$ in the form $(g_i)_{i\in I}$, where $g_i\in G_i$. If $I_1\subseteq I$, then for simplicity, we identify the subgroup $\{(g_i)_{i\in I}\in \prod_{i\in I}G_i\,|\,g_i=0$ for all $i\notin I_1\}$ of the group $\prod_{i\in I}G_i$ with the group $\prod_{i\in I_1}G_i$; we write elements of this group in the form $(g_i)_{i\in I_1}$.

As usual, $\mathbb{N}$ and $\mathbb{P}$ are the sets of positive integers and all prime integers, respectively, $\mathbb{Z}$ is the group (the ring) of integers, $\mathbb{Q}$ is the group (the field) of rational numbers. If $R$ is a unital ring, then $Re$ is the cyclic module over $R$ generated by the element $e$. If $S$ is a finite subset in $\mathbb{Z}$, then  $\text{gcd}(S)$ is the greatest common divisor of all integers in $S$ and $\text{lcm}(S)$ is the least common multiple of the integers in $S$. If $P_1\subseteq \mathbb{P}$, then \textsf{$P_1$-integer} is an integer such that any prime divisor of it (if it exists) is contained in $P_1$. It follows from the definition that $1$ is a $P_1$-integer for any $P_1\subseteq \mathbb{P}$. 
For any type $\tau$, we set
$$
P_{\infty}(\tau)=\{p\in \mathbb{P}\,|\, \tau(p)=\infty\}, \quad 
P_0(\tau)=\mathbb{P}\setminus P_{\infty}(\tau).
$$
Unless otherwise stated, we use notation and definitions from \cite{Fuc73}, \cite{Fuc15} and \cite{KMT03}.

\section{Multiplication Groups of Block-Rigid $CRQ$-groups of Ring Type}

\textsf{All over this section, $G$ is a reduced block-rigid $CRQ$-group of ring type with regulator $A$, regulator quotient $G/A=\langle d+A\rangle$ where $d\in G$, regulator index $n$ and set of critical types $T(G)=T(A)$.}

By setting $\text{Reg}_{\tau}\,G=A_{\tau}$, we can represent the group $A$ in the form $A=\oplus_{\tau\in T(G)}A_{\tau}$. According to \cite[Proposition 2.4.11]{Mad00}, such decomposition of the completely decomposable group $A$ is unique if and only if $A$ is a block-rigid group. For divisible hulls $\tilde G$, $\tilde A$, $\tilde A_{\tau}$ of the groups $G$, $A$ and $A_{\tau}$, respectively, we have relations 
$$
\tilde G=\tilde A=\oplus_{\tau\in T(G)}\tilde A_{\tau}.
$$
For $\tau\in T(G)$, we denote by $\pi_{\tau}$ the natural projection from the group $\tilde G$ onto $\tilde A_{\tau}$. 

In \cite{DugO93}, positive integers $m_{\tau}=m_{\tau}(G)$ ($\tau\in T(G)$) are defined; these integers are invariants of near-isomorphism of the group $G$. We can define integers $m_{\tau}$ ($\tau\in T(G)$ as follows; we take an element $d\in G/A$ such that $\langle d+A\rangle=G/A$. Let $d_{\tau}=\pi_{\tau}(d)\in \tilde{A}_{\tau}$ and let $m_{\tau}=o(d_{\tau}+A)$ be the order of the element $d_{\tau}+A$ in the torsion group $\tilde{A}/A$.
In \cite{DugO93}, it is shown that integers $m_{\tau}$ ($\tau\in T(G)$) do not depend on the choice of the element $d$. We note that $n=o(d+A)=\text{lcm}\{m_{\tau}\,|\,\tau\in T(G)\}$.

\textbf{Remark 2.1.}\\
Let $T$ be a finite set of pair-wise incomparable types and let $\{m_{\tau}\,|\,\tau\in T\}$ be some set of positive integers. We say that the set $\{m_{\tau}\,|\,\tau\in T\}$ satisfies \textsf{condition $(m)$} if for any $p\in \mathbb{P}$, $k\in \mathbb{N}$, $\tau\in T$, we have that $p^k$ divides $m_{\sigma}$ for some $\sigma\in T\setminus\{\tau\}$ provided $p^k$ divides $m_{\tau}$. We note that the set $\{m_{\tau}\,|\,\tau\in T\}$ satisfies condition $(m)$ if and only if the set $\{m_{\tau}\,|\,\tau\in T,\, m_{\tau}>1\}$ satisfies condition $(m)$.

According to \cite[Theorem 13.1.2]{Mad00}, the set $\{m_{\tau}\,|\,\tau\in T\}$ is a system of invariants of a near-isomorphism of some block-rigid $CRQ$-group $G$ with $T(G)=T$ if and only if this set satisfies condition $(m)$ and $m_{\tau}$ are $P_0(\tau)$-integers for all $\tau\in T$.~$\rhd$

In \cite[Theorem 3.5]{BlaM94}, it is proved that for any of the group $G\in \mathcal{A}_0$, there exists a direct decomposition 
$$
G=G_1\oplus C,\eqno (1)
$$
where $C$ is a completely decomposable group and $G_1$ is a rigid $CRQ$-group which satisfies the following conditions:
$$
\tau\in T(G_1) \text{ if and only if } m_{\tau}(G)> 1,\eqno (1')
$$
$$
m_{\tau}(G_1)=m_{\tau}(G) \text{ for all } \tau\in T(G_1).\eqno (1'')
$$
Decomposition $(1)$, which satisfies conditions $(1')$ and $(1'')$, is called a \textsf{main decomposition} of the group $G$. In a main decomposition of the group $G$, the group $G_1$ does not contain a completely decomposable summand; such groups are said to be \textsf{clipped}. We note that a main decomposition of a $CRQ$-group is not uniquely defined, since it depends on the choice of the element $d$ participating in the definition of the group. 
In what follows, we assume that a main decomposition of the group $G$ is fixed. We set $T_0(G)=\{\tau\in T(G)\,|\, m_{\tau}>1\}$. Then $T_0(G)$ is the set of critical types of a clipped direct summand in any main decomposition of the group $G$.

Let $B$ be the regulator of the group $G_1$, then $T(G_1)=T(B)=T_0(G)$ and $\tilde{G}_1=\tilde{B}$. There exists a system $E_0=\{e_0^{(\tau)}\in B_{\tau}\,|\, \tau\in T(B)\}$ such that
$$
B=\oplus_{\tau\in T(B)}R_{\tau}e_0^{(\tau)}. \eqno (2)
$$
In $(2)$, we assume that $R_{\tau}$ is a unitary subring in $\mathbb{Q}$, the type of the additive group of $R_{\tau}$ is equal to $\tau$, characteristics $\chi(e_0^{(\tau)})\in \tau$ contain only zeros and symbols $\infty$ ($\tau\in T(B)$). 

Let $D=\{d\in G_1\,|\,G/A=\langle d+A\rangle\}$, it is easy to see that $D\ne \varnothing$. Let $d\in D$. In the group $\tilde{B}$, the element $d$ can be represented in the form
$d=\sum_{\tau\in T(B)}\dfrac{s_{\tau}}{r_{\tau}}e_0^{(\tau)}$, where $s_{\tau}\in \mathbb{Z}$, $r_{\tau}\in \mathbb{N}$, $\text{gcd}(s_{\tau},r_{\tau})=1$.
Without loss of generality, we can assume that $s_{\tau}$, $r_{\tau}$ are $P_0(\tau)$-integers (otherwise, we can replace the system $E_0$). 

Let $\tau\in T(B)$. By the definition of the integer $m_{\tau}$, the relation $o\left(\dfrac{s_{\tau}}{r_{\tau}}e_0^{(\tau)}+A\right)=m_{\tau}$ holds in the group $\tilde{A}/A$. Since $r_{\tau}$ is a $P_0(\tau)$-integer and $\text{gcd}(s_{\tau},r_{\tau})=1$, we have $r_{\tau}=m_{\tau}$.
Consequently, the element $d$ of $\tilde{B}$ is of the form
$$
d=\sum_{\tau\in T(B)}\dfrac{s_{\tau}}{m_{\tau}}e_0^{(\tau)},
\eqno (3)
$$
and the integers $n$, $m_{\tau}$ and $s_{\tau}$ satisfy the following conditions:

$n=\text{lcm}\{m_{\tau}\,|\,\tau\in T(B)\}$,\hfill$(3')$

$\text{gcd}(s_{\tau},m_{\tau})=1$ for all $\tau\in T(B)$, \hfill$(3'')$

$s_{\tau}$ and $m_{\tau}$ are $P_0(\tau)$ numbers for any $\tau\in T(B)$.\hfill$(3''')$

A system $E_0=\{e_0^{(\tau)}\in B_{\tau} \,|\, \tau\in T(B)\}$ which satisfies conditions $(2)$ and $(3)$, is called a \textsf{$B$-basis} of the group $G$ defined by the element $d$. We note that the pair $(d,E_0)$ uniquely defines the numbers $s_{\tau}$ ($\tau\in T(B)$).
Relation $(3)$ is called a \textsf{standard representation} of block-rigid $CRQ$-group $G$ related to the pair $(d,E_0)$. 

\textbf{Remark 2.2.}\\
We note that a $B$-basis $E_0$ can be defined by more than one element $d\in D$.

$\lhd$ Indeed, let we have a standard representation $(3)$ of the group $G$. Let $\gamma$ be an integer which is co-prime with the regulator index $n$ and $d_1=\gamma d\in G_1$. Then $G/A=\langle d+A\rangle= \langle d_1+A\rangle$, i.e., $d_1\in D$. In addition,
$$
d_1=\sum_{\tau\in T(B)}\dfrac{\gamma s_{\tau}}{m_{\tau}}e_0^{(\tau)}. \eqno (4)
$$
If $\gamma$ is a $P_0(\tau)$-integer, then the relation $(4)$ is a standard representation of the group $G$. Consequently, the $B$-basis $E_0$ is defined by each of elements $d$ and $d_1$.

We note that if $\gamma$ is not a $P_0(\tau)$-integer, then the relation $(4)$ is not a standard representation of the group $G$.~$\rhd$

For a $B$-basis $E_0$, we set
$$
D(E_0)=\{d\in D\,|\, B\text{-basis } E_0 \text{ is determined by the element } d\}. 
$$ 
It follows from the definition of the $B$-basis that $D(E_0)\ne \varnothing$ and it follows from Remark 2.2 that $D(E_0)$ can contain more than one element.

With the right choice of elements $e_i^{(\tau)}\in C_{\tau}$ ($i=0,1,\ldots,k_{\tau}$), the group $C$ can be written in the form
$$
C=\oplus_{\tau\in T(C)}C_{\tau}=
\oplus_{\tau\in T(C)}\oplus_{i=1,\ldots,k_{\tau}}R_{\tau}e_i^{(\tau)},
$$
where $R_{\tau}$ is a unitary subring in the field of rational numbers, the type of the additive group of $R_{\tau}$ is equal to $\tau$, and characteristics $\chi(e_i^{(\tau)})\in \tau$ contain only zeros and symbols $\infty$. 

For $\tau\in T(G)$, we define the following sets:
$$
I_{\tau}(B)=
\begin{cases}
\{0\},\quad \text{ for } \tau\in T(B)\\
\quad \emptyset,\quad \text{ for } \tau\notin T(B),
\end{cases}
$$
$$
I_{\tau}(C)=
\begin{cases}
\{1,\dots,k_{\tau}\},\quad \text{ for } \tau\in T(C),\, k_{\tau}\in\mathbb{N}\\
\qquad \emptyset,\qquad \text{ for } \tau\notin T(C),
\end{cases}
$$
$$
I_{\tau}=I_{\tau}(B)\cup I_{\tau}(C).
$$
Then 
$A_{\tau}=\oplus_{i\in I_{\tau}}R_{\tau}e_i^{(\tau)}$ 
for any $\tau\in T(G)$ and 
$$
A=\oplus_{\tau\in T(G)}\oplus_{\tau\in I_{\tau}}R_{\tau}e_i^{(\tau)}. \eqno (5)
$$

A system $E=\{e_i^{(\tau)}\in A_{\tau} \,|\, \tau\in T(G), i\in I_{\tau}\}$ is called an \textsf{$A$-basis} of the group $G$ if $E$ satisfies $(5)$ and the subsystem $E_0=\{e_0^{(\tau)}\in A_{\tau} \,|\, \tau\in T(B)\}$ of $E$ is a $B$-basis.
 
Let $(G,\times)$ be a ring on the group $G\in\mathcal{A}_0$. Since $A$ is a fully invariant subgroup of the group $G$, we have that $A$ is an ideal of the ring $(G,\times)$ which is a direct sum of ideals $A_{\tau}$ ($\tau\in T(G)$). 
Thus, every multiplication on $G$ induces a multiplication on $A$; therefore, $\text{Mult}\,G
\subseteq\text{Mult}\,A$; however, the converse is not true. 

Let $E=\{e_i^{(\tau)}\in A_{\tau} \,|\, \tau\in T(G), i\in I_{\tau}\}$ be an $A$-basis of the group $G$. Then
for any set $\{u_{ij}^{(\tau)}\in A_{\tau}\,|\,\tau\in T(G),\;i,j\in I_{\tau}\}$, there exists a unique ring $(A,\times)$ such that $e_i^{(\tau)}\times e_j^{(\tau)}=u_{ij}^{(\tau)}$ for all $\tau\in T(G)$ and $i,j\in I_{\tau}$. 
The multiplication $\times$ is uniquely extended to a multiplication on $\tilde A=\tilde G$, where it is defined as follows:
$$
\sum_{i\in I_{\tau}}r_ie_i^{(\tau)}\times
\sum_{i\in I_{\tau}}r'_ie_i^{(\tau)}=\sum_{i,j\in I_{\tau}}r_ir'_j(e_i^{(\tau)}\times e_j^{(\tau)}) \eqno (6)
$$
for all $\tau\in T(G)$, $r_i,r'_j\in\mathbb{Q}$; and $\tilde A_{\tau}\times\tilde A_{\sigma}=0$ for $\tau\ne\sigma$.
However, $G$ is not necessarily a subring of the ring $(\tilde{A},\times)$. 
We say that the set $\{u_{ij}^{(\tau)}\in A_{\tau}\,|\,\tau\in T(G),\;i,j\in I_{\tau}\}$ \textsf{defines a multiplication} on $G$ with respect to the $A$-basis $E$ if there exists a ring $(G,\times)$ such that $e_i^{(\tau)}\times e_j^{(\tau)}=u_{ij}^{(\tau)}$ for all $\tau\in T(G)$ and $i,j\in I_{\tau}$. We note that any set $\{u_{ij}^{(\tau)}\in A_{\tau}\,|\,\tau\in T(G),\;i,j\in I_{\tau}\}$ defines multiplication on group $G$ with respect to the $A$-basis $E$ at most one way.

Let $G\in\mathcal{A}_0$. To describe sets defining multiplications on the group $G$, we define the following groups. 

For any $\tau\in T(G)$, let $n_{\tau}=|I_{\tau}|$ and let $M_{\tau}^{(0)}=M_{n_{\tau}}(A_{\tau})$ be the additive group of square matrices of order $n_{\tau}$ with elements in $A_{\tau}$,
$$
M_{\tau}^{(1)}=\begin{bmatrix}
m_{\tau}A_{\tau}&m_{\tau}A_{\tau}&\ldots&m_{\tau}A_{\tau}\\
m_{\tau}A_{\tau}&A_{\tau}&\ldots&A_{\tau}\\
\ldots&\ldots&\ldots&\ldots\\
m_{\tau}A_{\tau}&A_{\tau}&\ldots&A_{\tau}
\end{bmatrix},
$$
where the symbol $[\ldots]$ means the set of matrices of certain form, $m_{\tau}$ are invariants of near-isomorphism of the group $G$,
$$
M_{\tau}^{(2)}=\begin{bmatrix}
m_{\tau}^2A_{\tau}&m_{\tau}A_{\tau}&\ldots&m_{\tau}A_{\tau}\\
m_{\tau}A_{\tau}&A_{\tau}&\ldots&A_{\tau}\\
\ldots&\ldots&\ldots&\ldots\\
m_{\tau}A_{\tau}&A_{\tau}&\ldots&A_{\tau}
\end{bmatrix}
\subseteq M_{\tau}^{(1)}.
$$
We set 
$$
M^{(0)}=\prod_{\tau\in T(G)}M_{\tau}^{(0)},\;
M^{(1)}=\prod_{\tau\in T(G)}M_{\tau}^{(1)},\;
M^{(2)}=\prod_{\tau\in T(G)}M_{\tau}^{(2)}.
$$
Then $M^{(2)}\subseteq M^{(1)}\subseteq M^{(0)}$ and 
$M^{(2)}\cong M^{(1)}\cong M^{(0)}\cong \text{Mult}\,A$.

For the standard representation $(3)$ of the group $G$ related to the pair $(d,E_0)$ and for every $\tau\in T(G)$, we consider elements
$$
X^{(\tau)}=X^{(\tau)}(d,E_0)=\begin{pmatrix}
m_{\tau}s_{\tau}^{-1}e_0^{(\tau)}&0&\ldots&0\\
0&0&\ldots&0\\
\ldots&\ldots&\ldots&\ldots\\
0&0&\ldots&0
\end{pmatrix}\in M_{\tau}^{(1)},\;\text{if } \tau\in T(B),
$$
where $s_{\tau}^{-1}$ is an integer which is inverse to $s_{\tau}$ modulo $m_{\tau}$,
$$
X^{(\tau)}=\begin{pmatrix}
0&0&\ldots&0\\
0&0&\ldots&0\\
\ldots&\ldots&\ldots&\ldots\\
0&0&\ldots&0
\end{pmatrix}\in M_{\tau}^{(1)},\;\text{if } \tau\notin T(B).
$$
We set 
$$
X=X(d,E_0)=\left(X^{(\tau)}\right)_{\tau\in T(G)}=\left(X^{(\tau)}\right)_{\tau\in T(B)}\in M^{(1)}.
$$
$$
M(d,E_0)=\langle X,M^{(2)}\rangle\subseteq M^{(1)}.
$$
We note that integral solutions of the congruence $s_{\tau}x\equiv 1\,(\text{mod}\,m_{\tau})$ form the residue class modulo $m_{\tau}$. Therefore, the set $M(d,E_0)$ does not depend on the choice of the integers $s_{\tau}^{-1}$ in the definition of $X$.

We also note that if $\tau\in T(C)\setminus T(B)$, then $m_{\tau}=1$ by $(1')$. Therefore,
$$
M_{\tau}^{(2)}=M_{\tau}^{(1)}=M_{\tau}^{(0)}
$$
in this case.
For every set $U=\{u_{ij}^{(\tau)}\in A_{\tau}\,|\,\tau\in T(G),\;i,j\in I_{\tau}\}$ and every $\tau\in T(G)$, we consider the matrix 
$$
U^{(\tau)}=\begin{pmatrix}
u_{i_1,i_1}^{(\tau)}&u_{i_1,i_2}^{(\tau)}&\ldots&u_{i_1,i_{n_{\tau}}}^{(\tau)}\\
u_{i_2,i_1}^{(\tau)}&u_{i_2,i_2}^{(\tau)}&\ldots&u_{i_2,i_{n_{\tau}}}^{(\tau)}\\
\ldots&\ldots&\ldots&\ldots\\
u_{i_{n_{\tau}},i_1}^{(\tau)}&u_{i_{n_{\tau}},i_2}^{(\tau)}&\ldots&u_{i_{n_{\tau}},i_{n_{\tau}}}^{(\tau)}
\end{pmatrix}\; \in M_{\tau}^{(0)},
$$
where $i_k\in I_{\tau}$, $i_1<i_2<\ldots<i_{n_{\tau}}$. We set 
 $\overline{U}=\left(U^{(\tau)}\right)_{\tau\in T(G)}\in M^{(0)}$.

\textbf{Remark 2.3.} Let $G\in\mathcal{A}_0$, $G=G_1\oplus C$ be a main decomposition of the group $G$, and let $\text{Reg}\,G=A$, $\text{Reg}\,G_1=B$. In \cite{Kom09}, it is proved that for any multiplication $\times$ on group $G$, we have
$$
B_{\tau}\times A\subseteq m_{\tau}A_{\tau}\, \text{ and }\, 
A\times B_{\tau}\subseteq m_{\tau}A_{\tau}\, \text{ for all }\,\tau\in T(B).\;\rhd
$$

\textbf{Theorem 2.4.} Let $G$ be a block-rigid $CRQ$-group of ring type with $A$-basis $E$ containing an $B$-basis $E_0$. Let $U=\{u_{ij}^{(\tau)}\in A_{\tau}\,|\,\tau\in T(G),\;i,j\in I_{\tau}\}$. Then the following conditions are equivalent.
\begin{enumerate}
\item[\textbf{1)}]
$U$ defines a multiplication on $G$ with respect to the $A$-basis $E$.
\item[\textbf{2)}]
$\overline{U}\in M(d,E_0)$ for any $d\in D(E_0)$.
\item[\textbf{3)}]
$\overline{U}\in M(d,E_0)$ for some $d\in D(E_0)$.
\end{enumerate}

$\lhd$ 1)\,$\Rightarrow$\,2). Let the set $U=\{u_{ij}^{(\tau)}\in A_{\tau}\,|\,\tau\in T(G),\;i,j\in I_{\tau}\}$ induce the multiplication $\times$ on $G$ with respect to the $A$-basis $E=\{e_i^{(\tau)}\in A_{\tau}\,|\,\tau\in T(G),\;i\in I_{\tau}\}$. 
It follows from Remark 2.3 that $u_{0i}^{(\tau)},u_{i0}^{(\tau)}\in m_{\tau}A_{\tau}$ and
$u_{00}^{(\tau)}=m_{\tau}v_{00}^{(\tau)}$ ($v_{00}^{(\tau)}\in A_{\tau}$) for all $\tau\in T(B)$, $i\in I_{\tau}$.

Let $d\in D(E_0)$ and let a standard representation related to the pair $(d,E_0)$ be of the form
$$
d=\sum_{\tau\in T(B)} \dfrac{s_{\tau}}{m_{\tau}}e_{0}^{(\tau)}. 
$$
Since $d\times d\in G$, we have that $d\times d=\alpha d+a$ for some $\alpha\in\mathbb{Z}$ and $a\in A$. Then
$$
d\times d=\sum_{\tau\in T(B)}\dfrac{\alpha s_{\tau}}{m_{\tau}}e_0^{(\tau)}+ a. \eqno (7)
$$
On the other hand, it follows from $(6)$ that
$$
d\times d=\sum_{\tau\in T(B)}\dfrac{s_{\tau}}{m_{\tau}}e_{0}^{(\tau)}\times \sum_{\tau\in T(B)}\dfrac{s_{\tau}}{m_{\tau}}e_{0}^{(\tau)}=\sum_{\tau\in T(B)}\dfrac{s_{\tau}^2}{m_{\tau}^2}u_{00}^{(\tau)}= \sum_{\tau\in T(B)}\dfrac{s_{\tau}^2}{m_{\tau}}v_{00}^{(\tau)}. \eqno (8)
$$
Let $\tau\in T(B)$. It follows from $(7)$ and $(8)$ that
$$
\pi_{\tau}(d\times d)=\dfrac{s_{\tau}^2}{m_{\tau}}v_{00}^{(\tau)}=\dfrac{\alpha s_{\tau}}{m_{\tau}}e_{0}^{(\tau)}+a_{\tau},\;
\text{where}\;a_{\tau}=\pi_{\tau}(a)\in A_{\tau}.
$$
Consequently, $s_{\tau}^2v_{00}^{(\tau)}=\alpha s_{\tau}e_{0}^{(\tau)}+m_{\tau}a_{\tau}$; therefore,
$v_{00}^{(\tau)}=\alpha s_{\tau}^{-1}e_{0}^{(\tau)}+m_{\tau}a_{\tau}'$ for some $a_{\tau}'\in A_{\tau}$, where $s_{\tau}^{-1}$ is an integer which is inverse to $s_{\tau}$ modulo $m_{\tau}$. Therefore,
$$
u_{00}^{(\tau)}=m_{\tau}v_{00}^{(\tau)}=\alpha m_{\tau}s_{\tau}^{-1}e_{0}^{(\tau)}+m_{\tau}^2a_{\tau}'.
$$
Consequently, $\overline{U}\in\alpha X+M^{(2)}\subseteq M(d,E_0)$.

2)\,$\Rightarrow$\,3). The implication is directly verified.

3)\,$\Rightarrow$\,1). Let $\overline{U}\in M(d,E_0)$ for some $d\in D(E_0)$ and let a standard representation related to $(d,E_0)$ be of the form $d=\sum_{\tau\in T(B)}\dfrac{s_{\tau}}{m_{\tau}}e_0^{(\tau)}$. Then $\overline{U}=\alpha X+Y$ for some $\alpha\in\mathbb{Z}$, $Y\in M^{(2)}$. Therefore, for all $\tau\in T(B)$ and $i\in I_{\tau}$, we have
$$
u_{0i}^{(\tau)}=m_{\tau}v_{0i}^{(\tau)},\quad u_{i0}^{(\tau)}=m_{\tau}v_{i0}^{(\tau)},\quad
u_{00}^{(\tau)}=\alpha m_{\tau}s_{\tau}^{-1}e_{0}^{(\tau)}+m_{\tau}^2a_{\tau},
$$
where $v_{0i}^{(\tau)}, v_{i0}^{(\tau)}, a_{\tau}\in A_{\tau}$ and the integers $s_{\tau}^{-1}$ satisfy conditions $s_{\tau}s_{\tau}^{-1}=1+m_{\tau}x_{\tau}$ for some $x_{\tau}\in \mathbb{Z}$.

There exists a ring $(A,\times)$ such that $e_i^{(\tau)}\times e_j^{(\tau)}=u_{ij}^{(\tau)}$ for all $\tau\in T(G)$, $i,j\in I_{\tau}$; and $A_{\tau}\times A_{\sigma}=0$ for $\tau\ne\sigma$. 
This multiplication is extended to a multiplication on divisible hull $\tilde A=\tilde G$ of the group $A$. We prove that $G$ is a subring of the ring $(\tilde A,\times)$. 

It follows from $(6)$ that
$$
d\times d=\sum_{\tau\in T(B)}\dfrac{s_{\tau}^2}{m_{\tau}^2}u_{00}^{(\tau)}=\alpha\sum_{\tau\in T(B)}\dfrac{s_{\tau}^2s_{\tau}^{-1}}{m_{\tau}}e_0^{(\tau)}+ \sum_{\tau\in T(B)}s_{\tau}^2a_{\tau}=
$$
$$
=\alpha\sum_{\tau\in T(B)}\dfrac{s_{\tau}}{m_{\tau}}(1+m_{\tau}x_{\tau})e_0^{(\tau)}+ \sum_{\tau\in T(B)}s_{\tau}^2a_{\tau}=
$$
$$
=\alpha\sum_{\tau\in T(B)}\dfrac{s_{\tau}}{m_{\tau}}e_0^{(\tau)} +\sum_{\tau\in T(B)}(\alpha s_{\tau}x_{\tau})e_0^{(\tau)}+ \sum_{\tau\in T(B)}s_{\tau}^2a_{\tau}=
$$
$$
=\alpha d+ \sum_{\tau\in T(B)}(\alpha s_{\tau}x_{\tau}e_0^{(\tau)}+s_{\tau}^2a_{\tau})\in G.
$$
In addition, if $\sigma\in T(B)$ and $i\in I_{\sigma}$, then
$$
d\times e_i^{(\sigma)}=\left(\sum_{\tau\in T(B)}\dfrac{s_{\tau}}{m_{\tau}}e_0^{(\tau)}\right)\times e_i^{(\sigma)}=\dfrac{s_{\sigma}}{m_{\sigma}}u_{0i}^{(\sigma)}=
\dfrac{s_{\sigma}}{m_{\sigma}}m_{\sigma}v_{0i}^{(\sigma)}=s_{\sigma}v_{0i}^{(\sigma)}\in A.
$$
Similarly, we have $e_i^{(\sigma)}\times d\in A$.

If $\sigma\notin T(B)$, then $d\times e_i^{(\sigma)}= e_i^{(\sigma)}\times d=0$ for any $i\in I_{\sigma}$. Since $G=\langle d,A\rangle$, we have that $G$ is a subring of the ring $(\tilde A,\times)$. Therefore, the set $U$ defines a multiplication on $G$.~$\rhd$

It follows from Theorem 2.4 that the group $M(d,E_0)$ does not depend on the choice of the element $d\in D(E_0)$. In the following assertion, we consider relations between elements of groups $M(d_1,E_0)$ and $M(d_2,E_0)$ for $d_1,d_2\in D(E_0)$. We note that if $d_1,d_2\in D$, then $\langle d_1+A\rangle=\langle d_2+A\rangle=G/A$; therefore, $d_1=\gamma d_2+b$ for  some $b\in B$ and some integer $\gamma$ which is co-prime with $n(G)$.

\textbf{Proposition 2.5.} Let $G$ be a group in $\mathcal{A}_0$ with main decomposition $(1)$ and with a $B$-basis $E_0$ and regulator index $n$. Then the following assertions hold.
\begin{enumerate}
\item[\textbf{1)}]
$M(d_1,E_0)=M(d_2,E_0)$ for any $d_1,d_2\in D(E_0)$.
\item[\textbf{2)}]
If $d_1,d_2\in D(E_0)$ and $d_1=\gamma d_2+b$, where $\gamma\in \mathbb{Z}$, $\text{lcd}(\gamma,n)=1$, $b\in B$, $X_1=X(d_1,E_0)\in M(d_1,E_0)$, $X_2=X(d_2,E_0)\in M(d_2,E_0)$, then 
$X_1+M^{(2)}=\gamma^{-1}X_2+M^{(2)}$, where $\gamma^{-1}$ is an integer which is inverse to $\gamma$ modulo $n$.
\end{enumerate}

$\lhd$. Let $d_1,d_2\in D(E_0)$, $d_1=\gamma d_2+b$, where $\gamma\in\mathbb{Z}$, $\text{lcd}\,(\gamma,n)=1$ and $b=\sum_{\tau\in T(B)}b_{\tau}e_0^{(\tau)}$ ($b_{\tau}\in R_{\tau}$). Let
$$
d_1=\sum_{\tau\in T(B)}\dfrac{s_{\tau}}{m_{\tau}}e_0^{(\tau)},\quad 
d_2=\sum_{\tau\in T(B)}\dfrac{t_{\tau}}{m_{\tau}}e_0^{(\tau)}
$$
be standard representations of the group $G$ related to $(d_1,E_0)$ and $(d_2,E_0)$, respectively. In the divisible hull $\tilde G$ of the group $G$, we have relations
$$
\sum_{\tau\in T(B)}\dfrac{s_{\tau}}{m_{\tau}}e_0^{(\tau)}=d_1=\gamma d_2+b=\sum_{\tau\in T(B)}\dfrac{\gamma t_{\tau}+m_{\tau}b_{\tau}}{m_{\tau}}e_0^{(\tau)}.
$$
Therefore,
$$
s_{\tau}=\gamma t_{\tau}+m_{\tau}b_{\tau}\;\text{for all}\; \tau\in T(B).\eqno(9)
$$
Let $\tau\in T(B)$, $\gamma^{-1}$ be an integer which is inverse to $\gamma$ modulo $n$, and let $t_{\tau}^{-1}$ be an integer which is inverse to $t_{\tau}$ modulo $m_{\tau}$. Then number $\gamma^{-1}$ is inverse to $\gamma$ modulo $m_{\tau}$ by $(3')$; therefore, the integer $\gamma^{-1}t_{\tau}^{-1}$ is inverse to $s_{\tau}$ modulo $m_{\tau}$ by $(9)$.


Let $a_{\tau}\in A_{\tau}$. Then
$$
s_{\tau}^{-1}m_{\tau}e_0^{(\tau)}+m_{\tau}^2a_{\tau}=
\gamma^{-1}t_{\tau}^{-1}m_{\tau}e_0^{(\tau)}+m_{\tau}^2a_{\tau}',\;\text{where}\; a_{\tau}'\in A_{\tau}.
$$
Consequently, 
$$
X_1+M^{(2)}\subseteq \gamma^{-1}X_2+M^{(2)}.\eqno (10)
$$
Since $d_2=\gamma^{-1}d_1+b_1'$ (where $b_1'\in B$), we have $$
\gamma^{-1}X_2+M^{(2)}\subseteq \gamma^{-1}(\gamma X_1+M^{(2)})+M^{(2)}=X_1+M^{(2)}
$$
by $(10)$.~$\rhd$

It follows from Proposition 2.5 that we can write $M(d,E_0)=M(E_0)$, however we will often write $M(d,E_0)$ if we want to point which a standard representation is used in the definition of the group $M(d,E_0)$.

\textbf{Remark 2.6.}\\
Let $\overline{U}\in M^{(2)}$ for the set $U=\{u_{ij}^{(\tau)}\in A_{\tau}\,|\,\tau\in T(G),\;i,j\in I_{\tau}\}$. It follows from Theorem 2.4 that, with respect to any $A$-basis $E$ of the group $G$, the set $U$ defines a multiplication $\times_{U,E}$ on $G$ such that $G\times_{U,E}G\subseteq A$. Such a multiplication is called a \textsf{regulator} multiplication.~$\rhd$

\textbf{Example 2.7.} It is possible that the set $U=\{u_{ij}^{(\tau)}\,|\,\tau\in T(G),\;i,j\in I_{\tau}\}$ defines a multiplication on $G$ with respect to one $A$-basis and does not define any multiplication on $G$ with respect to another $A$-basis even for the same main decomposition. Moreover, the following situation is possible: there exist two $A$-bases $E$ and $F$ such that any set $U$, which defines non-regulator multiplication with respect to $E$, does not define any multiplication with respect to $F$. This means that
$$
M(E_0)\cap M(F_0)=M^{(2)}.
$$
$\lhd$ Let $s_1$, $s_2$ be two co-prime integers, $s_1>1$, $s_2>1$, and let $m$ be a prime integer which does not divide any of the integers $s_1$, $s_2$, $s_1^2-s_2^2$.

Let $\tau_i$ be an idempotent type such that $P_0(\tau_i)$ is the set of all prime divisors of the integers $m$ and $s_i$, respectively ($i=1,2$). Then types $\tau_1$ and $\tau_2$ are incomparable.

We consider a group $B=R_1e_1\oplus R_2e_2$, where $R_1$ and $R_2$ are unital subrings in the field $\mathbb{Q}$ whose additive groups are of types $\tau_1$ and $\tau_2$, respectively.

It follows from Remark 2.1 that there exists a $CRQ$-group $G=\langle d,B\rangle$ with regulator $B$ and quasi-isomorphism invariants $m_{\tau_1}=m_{\tau_2}=m$. We can choose the group $G$ in such a way that a standard representation of $G$ is of the form
$$
d=\dfrac{s_1}{m}e_1+\dfrac{s_2}{m}e_2.
$$
Then the system $E_0=\{e_1,e_2\}$ is a $B$-basis of the group $G$ defined by the element $d$ (in this case, the $B$-basis coincides with an $A$-basis).

We set
$$
B_1=B_{\tau_1}=R_1e_1,\; B_2=B_{\tau_2}=R_2e_2.
$$
We consider $d_1=d+(e_1+e_2)\in G$. Then $d_1\in D$;  in $\tilde G$, we have
$$
d_1=\dfrac{s_1+m}{m}e_1+\dfrac{s_2+m}{m}e_2.\eqno (11)
$$
Since $s_i+m$ is co-prime with each of the integers $s_i$ and $m$, we have that $s_i+m$ is a $P_{\infty}(\tau_i)$-integer for $i=1,2$. Consequently, $(11)$ is not a standard representation of the group $G$. We set $f_1=(s_1+m)e_1$ and $f_2=(s_2+m)e_2$. Then the system $F_0=\{f_1,f_2\}$ is a $B$-basis defined by the element $d_1$. The standard representation of the group $G$, related to the pair $(d_1,F_0)$, is of the form $d_1=\dfrac{1}{m}f_1+\dfrac{1}{m}f_2$.

Let the set $U=\{u_{i}\in B_i\,|\,i=1,2\}$ define a non-regulator multiplication on $G$ with respect to of the $B$-basis $E_0=\{e_1,e_2\}$. By Theorem 2.4, we have $(u_1,u_2)\in M(d,E_0)\setminus M^{(2)}$.
Consequently,
$$
u_1\in \alpha ms_1^{-1}e_1+m^2B_1,\quad 
u_2\in \alpha ms_2^{-1}e_2+m^2B_2,\eqno (12)
$$
for some integer $\alpha$ which is not divided by the prime integer $m$. 

We assume that the set $U$ defines multiplication with respect to the $B$-basis $F_0=\{f_1,f_2\}$. Then it follows from Theorem 2.4 that for some $\beta\in\mathbb{Z}$, we have
$$
u_1\in \beta mf_1+m^2B_1,\quad 
u_2\in \beta mf_2+m^2B_2.\eqno (13)
$$
It follows from $(12)$ and $(13)$ that
$$
\alpha ms_1^{-1}e_1\in\beta mf_1+m^2B_1=\beta ms_1e_1+m^2B_1,
$$
$$
\alpha ms_2^{-1}e_2\in\beta mf_2+m^2B_2=\beta ms_2e_2+m^2B_2.
$$
Therefore,
$$
\alpha= \beta s_1^2+ mx_1, \quad \alpha= \beta s_2^2+ mx_2 \eqno (14)
$$
for some $x_1\in R_1$, $x_2\in R_2$. Since $x_i=\dfrac{\alpha-\beta s_i^2}{m}\in R_i$ and $m$ is a $P_0(\tau_i)$-integer, we have $x_i\in\mathbb{Z}$ for $i=1,2$. Consequently, it follows from $(14)$ that $\beta (s_1^2-s_2^2)=my$ for $y=x_2-x_1\in\mathbb{Z}$. Since the prime integer $m$ does not divide $s_1^2-s_2^2$, we have that $m$ divides $\beta$; therefore, $m$ divides $\alpha$ by $(14)$. This contradicts to the property that $(u_1,u_2)\notin M^{(2)}$. Consequently, the set $U=\{u_1,u_2\}$ does not define any multiplication with respect to the $B$-basis $F_0$.~$\rhd$

Let $\times$ be a multiplication on a group $G\in \mathcal{A}_0$. Let $E=\{e_i^{(\tau)}\in A_{\tau}\,|\,\tau\in T(G),\,i\in I_{\tau}\}$ be an $A$-basis of the group $G$ and let
$$
U_{\times}=U_{\times}(E)=\{u_{ij}^{(\tau)}=e_{i}^{(\tau)}\times e_{j}^{(\tau)}\in A_{\tau}\,|\,\tau\in T(G),\;i,j\in I_{\tau}\}.
$$
It clearly follows from Theorem 2.4 that the correspondence $\times\mapsto \overline{U_{\times}}$ defines an isomorphism from the group $\text{Mult}\,G$ onto $M(E_0)$.

\textbf{Theorem 2.8.} If $G\in \mathcal{A}_0$ and $E_0$ is a $B$-basis of the group $G$, then $\text{Mult}\,G\cong M(E_0)$.~$\rhd$

We note that Theorem 2.8 implies the following property: up to isomorphism, the group $M(E_0)$ does not depend on the choice of the $B$-basis $E_0$.

\textbf{Remark 2.9.}\\
Let $G=\langle d,A\rangle\in \mathcal{A}_0$ and let $E$ be an $A$-basis of the group $G$ containing the $B$-basis $E_0$.

\textbf{1.} It follows from Theorem 2.8 that the group $\text{Mult}\,G$ can be identified with the group $M(E_0)=\langle X,M^{(2)}\rangle$ and the multiplication $\times$ can be identified with $\overline{U}_{\times}\in M(E_0)$.

\textbf{2.} Let $\overline{U}_{\times}=\overline{U}_{\times}(E)\in\text{Mult}\,G=M(d,E_0)=\langle ,M^{(2)}\rangle$, $\alpha\in\mathbb{Z}$. It follows from the proof of Theorem 2.4 that $\overline{U}_{\times}\in \alpha X+M^{(2)}$ if and only if $d\times d\in \alpha d+A$.

\textbf{3.} It follows from \textbf{2} that $\overline{U}_{\times}\in M^{(2)}$ if and only if $G\times G\subseteq A$. In the group $\text{Mult}\,G$, this means that the subgroup $\text{Hom}(G\otimes G,A)$ of all regulator multiplications coincides with the group $M^{(2)}$.~$\rhd$

\section{Properties of Multiplication Groups of Block-Rigid $CRQ$-Groups of Ring Type}\label{section3}

The purpose of this section is to show that for any group $G$ in the class $\mathcal{A}_0$, the group $\text{Mult}\,G$ belongs to this class, as well. We will also describe the rank, the set of critical types, invariants of near-isomorphism, the regulator, a main decomposition, and a standard representation of the group $\text{Mult}\,G$, where $G\in \mathcal{A}_0$.

\textbf{Remark 3.1.}\\ 
Let $A$ be a completely decomposable block-rigid group of finite rank and $G=\langle d,A\rangle$, where $d\in \tilde A$. Let we have  $o(d_{\tau}+A)=m_{\tau}$ in the group $\tilde A/A$, where $d_{\tau}=\pi_{\tau}(d)$ for $\tau\in T(A)$. 
Then the set $\{m_{\tau}\,|\,\tau\in T(A)\}$ satisfies condition $(m)$ (see Remark 2.1) if and only if for any $\tau\in T(A)$, the subgroup $A_{\tau}$ is pure in $G$.

$\lhd$
Indeed, let the set $\{m_{\tau}\,|\,\tau\in T(A)\}$ satisfy condition $(m)$,~ $\sigma\in T(A)$, $a\in A_{\sigma}$ and $a=k(td+x)$ for some $k,t\in \mathbb{Z}$ and $x\in A$. Let $\tau\ne \sigma$, then
$$
ktd_{\tau}+kx_{\tau}=0,\; \text{ where }\,x_{\tau}=\pi_{\tau}(x)\in A_{\tau}.
$$
Therefore, $td_{\tau}\in A_{\tau}$, whence $m_{\tau}$ divides $t$. Since the set $\{m_{\tau}\,|\,\tau\in T(A)\}$ satisfies condition $(m)$, we have that
$$
\text{lcm}\{m_{\tau}\,|\,\tau\in T(A), \tau\ne\sigma\}=
\text{lcm}\{m_{\tau}\,|\,\tau\in T(A)\}=n.
$$
Consequently, $n$ divides $t$; therefore, $td+x\in A$. Since the type of the element $td+x$ is equal to $\sigma$, we have that $td+x\in A_{\sigma}$.

Conversely, let the subgroup $A_{\tau}$ be pure in $G$ for any $\tau\in T(A)$. We assume that the set $\{m_{\tau}\,|\,\tau\in T(A)\}$ does not satisfy condition $(m)$. Then there exists a type $\sigma\in T(A)$ such that $m_{\sigma}$ does not divide $n_1=\text{lcm}\{m_{\tau}\,|\,\tau\in T(A), \tau\ne\sigma\}$. Consequently, for $n=\text{lcm}\{m_{\tau}\,|\,\tau\in T(A)\}$, it is true that $n=n_1n_2$ for some integer $n_2>1$. Consequently, $n_1d_{\sigma}\notin A_{\sigma}$ and $n_2(n_1d_{\sigma})\in A_{\sigma}$. Therefore, the subgroup $A_{\sigma}$ is not pure in $G$.~$\rhd$

\textbf{Remark 3.2.}\\ 
Let $A$ be a reduced block-rigid (resp., rigid) completely decomposable group of finite rank and let $G=\langle d,A\rangle$, where $d\in\tilde A\setminus A$. Then $A$ is a subgroup of finite index of the group $G$, and $T(G)=T(A)$. Therefore, $G$ is a block-rigid (resp., rigid) $ACD$-group by the definition. We note that if $A$ is a group of ring type, then $G$ also is a group of ring type. Since $G/A=\langle d+A\rangle$ is a cyclic group, we have that $G$ is a $CRQ$-group by \cite[Section 2]{BlaM94}. In addition, $A=\text{Reg}\,G$ if and only if subgroup $A_{\tau}$ is pure in $G$ \cite[Section 2]{BlaM94} for any $\tau\in T(G)$.~$\rhd$

Let $G\in \mathcal{A}_0$. The following theorem describes properties of the group $\text{Mult}\,G$. In what follows,
$G\in \mathcal{A}_0$, $T(G)=T$, $T_0(G)=T_0$, $m_{\tau}(G)=m_{\tau}$ $(\tau\in T)$, $n(G)=n$, $G=G_1\oplus C$ is a main decomposition of the group $G$, $\text{Reg}\,G_1=B=\oplus_{\tau\in T_0}R_{\tau}e_0^{(\tau)}$, $\text{Reg}\,G=A=B\oplus C$, $G=\langle d,A\rangle$,
and a standard representation of the group $G$ is of the form
$$
d=\sum_{\tau\in T_0}\dfrac{s_{\tau}}{m_{\tau}}e_0^{(\tau)}.
$$
We note that the set of integral solutions of the congruence $s_{\tau}x\equiv 1 (\text{mod }m_{\tau})$ always contains a $P_0(\tau)$-integer $s_{\tau}^{\varphi(m_{\tau})-1}$, where $\varphi(x)$ is the Euler function. Therefore, we always can take this $P_0(\tau)$-integer as the integer $s_{\tau}^{-1}$ inverse to $s_{\tau}$ modulo $m_{\tau}$, . 

\textbf{Theorem 3.3.} Let $G\in \mathcal{A}_0$. Then the group $\text{Mult}\,G$ satisfies the following conditions.

\textbf{1.} The group $\text{Mult}\,G$ is a block-rigid $CRQ$-group of ring type with regulator $M^{(2)}=\text{Hom}(G\otimes G,A)$.

\textbf{2.} $T(\text{Mult}\,G)=T(G)$ and $T_0(\text{Mult}\,G)=T_0(G)$, as a corollary.

\textbf{3.} $m_{\tau}(\text{Mult}\,G)=m_{\tau}(G)$ for any $\tau\in T(G)$,~ $n(\text{Mult}\,G)=n(G)$.

\textbf{4.} $r(\text{Reg}_{\tau}(\text{Mult}\,G))=(r(\text{Reg}_{\tau}G))^3$ for any $\tau\in T(G)$.

\textbf{5.} One of main decompositions of the group $\text{Mult}\,G$ is of the form $\text{Mult}\,G=M'\oplus M''$, where
$$
M'=\langle X,K\rangle,\; K=\prod_{\tau\in T_0(G)}K_{\tau},\quad
K_{\tau}=\begin{bmatrix}
m_{\tau}^2B_{\tau}&0&\ldots&0\\
0&0&\ldots&0\\
\ldots&\ldots&\ldots&\ldots\\
0&0&\ldots&0
\end{bmatrix} \subseteq M_{\tau}^{(2)},
$$
$$
X=\left(X^{(\tau)}\right)_{\tau\in T_0(G)},\qquad
X^{(\tau)}=\begin{pmatrix}
m_{\tau}s_{\tau}^{-1}e_0^{(\tau)}&0&\ldots&0\\
0&0&\ldots&0\\
\ldots&\ldots&\ldots&\ldots\\
0&0&\ldots&0
\end{pmatrix}\in M_{\tau}^{(1)} \text{ for } \tau\in T_0(G),
$$
$$
M''=\prod_{\tau\in T(G)}M_{\tau}'',\qquad
M_{\tau}''=\begin{bmatrix}
m_{\tau}^2C_{\tau}&m_{\tau}A_{\tau}&\ldots&m_{\tau}A_{\tau}\\
m_{\tau}A_{\tau}&A_{\tau}&\ldots&A_{\tau}\\
\ldots&\ldots&\ldots&\ldots\\
m_{\tau}A_{\tau}&A_{\tau}&\ldots&A_{\tau}
\end{bmatrix}\, \subseteq M_{\tau}^{(2)}.
$$
In addition, $T(M')=T_0(G)$, $\text{Reg}\,M'=K$.

\textbf{6.} For every $\tau\in T_0(G)$, we denote by
$s_{\tau}^{-1}$ a $P_0(\tau)$-integer which is inverse to $s_{\tau}$ modulo $m_{\tau}$,
$$
E_0^{(\tau)}=\begin{pmatrix}
m_{\tau}^2e_0^{(\tau)}&0&\ldots&0\\
0&0&\ldots&0\\
\ldots&\ldots&\ldots&\ldots\\
0&0&\ldots&0
\end{pmatrix}\in K_{\tau}.
$$
Then the system $\{E_0^{(\tau)}\,|\,\tau\in T_0(G)\}$ is one of $B$-bases of the group $\text{Mult}\,G$. One of standard representations of the group $\text{Mult}\,G$ is of the form
$$
X=\left(\dfrac{s_{\tau}^{-1}}{m_{\tau}}E_0^{(\tau)}\right)_{\tau\in T_0(G)}.
$$

$\lhd$ \textbf{1.} It follows from Theorem 2.8 that $\text{Mult}\,G=\langle X,M^{(2)}\rangle$, where
$$
X=\left(X^{(\tau)}\right)_{\tau\in T_0},\;
X^{(\tau)}=\begin{pmatrix}
m_{\tau}s_{\tau}^{-1}e_0^{(\tau)}&0&\ldots&0\\
0&0&\ldots&0\\
\ldots&\ldots&\ldots&\ldots\\
0&0&\ldots&0
\end{pmatrix}\in M^{(1)}\; \text{ for } \tau\in T_0.
$$ 
Since $\text{gcd}(s_{\tau}^{-1},m_{\tau})=1$ for $\tau\in T_0$,  we have $o(X_{\tau}+M^{(2)})=m_{\tau}$ in the group $\tilde M^{(2)}/M^{(2)}$.

By Remark 2.1, the set $\{m_{\tau}\,|\,\tau\in T\}$ satisfies condition $(m)$. Therefore, it follows from Remark 3.1 that the subgroups $M_{\tau}^{(2)}$ are pure in $\text{Mult}\,G$ for any $\tau\in T$. Since $M^{(2)}$ is a block-rigid completely decomposable group of ring type, $\text{Mult}\,G$ is a block-rigid $CRQ$-group of ring type with regulator $M^{(2)}$ by Remark 3.2. It follows from Remark 2.9(3) that we have $M^{(2)}=\text{Hom}(G\otimes G,A)$.

\textbf{2.} It follows from \textbf{1} and the definition of the group $M^{(2)}$ that
$$
T(\text{Mult}\,G)=T(M^{(2)})=T(G).
$$
\textbf{3.} We have $\text{Mult}\,G=\langle X,M^{(2)}\rangle$ and $\text{Reg}(\text{Mult}\,G)=M^{(2)}$ by \textbf{1}. Let $\tau\in T$, then
$$
m_{\tau}(\text{Mult}\,G)=o(X_{\tau}+M^{(2)})=m_{\tau}=m_{\tau}(G).
$$
Therefore, $n(\text{Mult}\,G)=\text{lcm}\{m_{\tau}\,|\,\tau\in T\}=n(G)$.

\textbf{4.} Let $\tau\in T$. It follows from \textbf{1} that
$$
\text{Reg}_{\tau}(\text{Mult}\,G)=M_{\tau}^{(2)}=\begin{bmatrix}
m_{\tau}^2A_{\tau}&m_{\tau}A_{\tau}&\ldots&m_{\tau}A_{\tau}\\
m_{\tau}A_{\tau}&A_{\tau}&\ldots&A_{\tau}\\
\ldots&\ldots&\ldots&\ldots\\
m_{\tau}A_{\tau}&A_{\tau}&\ldots&A_{\tau}
\end{bmatrix}\cong M_{n_{\tau}}(A_{\tau}),
$$ 
where $n_{\tau}=r(A_{\tau})$. Consequently, $r(\text{Reg}_{\tau}(\text{Mult}\,G))=n_{\tau}^3=(r(\text{Reg}_{\tau}\,G))^3$.

\textbf{5.} In the decomposition $\text{Mult}\,G\cong M'\oplus M''$, the group $M''$ is completely decomposable and $M'=\langle X,K\rangle$. By the definition of the group $K$, we have $T(K)=T(M')=T_0$. It is easy to see that $o(X^{(\tau)}+K)=m_{\tau}$ in the group $\tilde K/K$, for any $\tau\in T_0$. 
Since $\{m_{\tau}\,|\,\tau\in T_0\}$ satisfies condition $(m)$ (by Remark 2.1) and $K$ is a rigid completely decomposable group, we have that $M'$ is a rigid group in $\mathcal{A}_0$ with $\text{Reg}\,M'=K$ by Remark 3.1 and Remark 3.2.

Since $T(M')=T_0$, we have that $\tau\in T(M')$ if and only if $m_{\tau}(\text{Mult}\,G)=m_{\tau}>1$. In addition,
$$
m_{\tau}(M')=o(X^{(\tau)}+K)=m_{\tau}=m_{\tau}(\text{Mult}\,G),
$$
by \textbf{3}. It follows from $(1')$, $(1'')$ that the decomposition $\text{Mult}\,G= M'\oplus M''$ is a main decomposition of the group $\text{Mult}\,G$.

\textbf{6.} Let $\tau\in T_0$, $s_{\tau}^{-1}$ be a $P_0(\tau)$-integer which is inverse to $s_{\tau}$ modulo $m_{\tau}$, and let
$$
E_0^{(\tau)}=\begin{pmatrix}
m_{\tau}^2e_0^{(\tau)}&0&\ldots&0\\
0&0&\ldots&0\\
\ldots&\ldots&\ldots&\ldots\\
0&0&\ldots&0
\end{pmatrix}\in K_{\tau}.
$$
Then $K=\prod_{\tau\in T_0}R_{\tau}E_0^{(\tau)}$. In addition, 
$$
X=\left(\dfrac{s_{\tau}^{-1}}{m_{\tau}}E_0^{(\tau)}\right)_{\tau\in T_0}.\eqno (15)
$$
Since $s_{\tau}^{-1}$ ($\tau\in T_0$) is a $P_0(\tau)$-integer, $(15)$ is a standard representation of the group $\text{Mult}\,G$. Therefore, $\{E_0^{(\tau)}\,|\,\tau\in T_0(G)\}$ is a $B$-basis of the group $\text{Mult}\,G$.~$\rhd$

\section{Data Availability Statement}

Our manuscript has no associate data.

\end{document}